\documentclass[12pt]{article}
\usepackage{amsmath,amssymb,amsbsy,amsfonts,amsthm,latexsym,
                        amsopn,amstext,amsxtra,euscript,amscd}

\begin{document}

\newtheorem{theorem}{Theorem}
\newtheorem{lemma}[theorem]{Lemma}
\newtheorem{claim}[theorem]{Claim}
\newtheorem{cor}[theorem]{Corollary}
\newtheorem{prop}[theorem]{Proposition}
\newtheorem{definition}{Definition}
\newtheorem{question}[theorem]{Open Question}

\baselineskip 15pt

\title{On a symmetric congruence and its applications}

\author{{M.~Z.~Garaev}
\\
\normalsize{Instituto de Matem{\'a}ticas,  UNAM}
\\
\normalsize{Campus Morelia, Apartado Postal 61-3 (Xangari)}
\\
\normalsize{C.P. 58089, Morelia, Michoac{\'a}n, M{\'e}xico} \\
\normalsize{\tt garaev@matmor.unam.mx} \\
\and\\
{A.~A.~Karatsuba}
\\
\normalsize{Steklov Institute of Mathematics}
\\
\normalsize{Russian Academy of Sciences}
\\
\normalsize{GSP-1, ul. Gubkina 8}
\\
\normalsize{Moscow, Russia} \\
\normalsize{\tt karatsuba@mi.ras.ru}
 }

\date{\empty}

\pagenumbering{arabic}

\maketitle


\begin{abstract}
For a large integer $m,$ we obtain an asymptotic formula for the
number of solutions of a certain congruence modulo $m$ with four
variables, where the variables belong to special sets of residue
classes modulo $m.$ This formula are applied to obtain a new bound
for a double trigonometric sum with an exponential function and new
information on the exceptional set of the multiplication table
problem in a residue ring modulo $m.$
\end{abstract}


\section{Introduction}
Throughout the paper the implied constants in the Landau `$O$' and
`$o$' symbols as well as in the Vinogradov symbols `$\ll$' and
`$\gg$' may depend on the small positive quantity $\varepsilon.$ By
$p$ and $q$ we will always denote prime numbers.

Let $m$ be an integer parameter, $\mathcal{V}$ be any subset of
prime numbers coprime to $m$ and not exceeding $m^{1/2},$ and let
$S$ and $L$ be any integers with $0<L\le m.$ In this paper we obtain
an asymptotic formula for the number of solutions of the congruence
$$
v_1y_1\equiv v_2y_2\pmod m,\quad v_1, v_2\in \mathcal{V},\quad
S+1\le y_1, y_2\le S+L
$$
and give its applications.

Below $|\mathcal{V}|$ stands for the number of elements in
$\mathcal{V}.$ By $\tau(m)$ we denote the classical divisor
function.

\begin{theorem}
\label{thm:main} {\it The following asymptotic formula holds:
$$
J=\frac{|\mathcal{V}|^2L^2}{m}+|\mathcal{V}|L-\frac{|\mathcal{V}|L^2}{m}+
O\left(\frac{m^2\log^2m}{\phi(m)}\right),
$$}
where $\phi(m)$ is the Euler function.
\end{theorem}

Theorem~\ref{thm:main} finds its application in estimating of double
trigonometric sums with an exponential function. For a positive
integer $m$ we denote by $\mathbb{Z}_m=\{0, 1,..., m-1\}$ the
residue ring modulo $m.$ Let $p$ be a large prime number, $T$ be a
divisor of $p-1,$  $\lambda$ be an element of $\mathbb{Z}_p$ of
multiplicative order $T,$ i.e. $\lambda=g^{(p-1)/T}$ for some
primitive root $g$ modulo $p,$ $\gamma(n)$ be any complex
coefficients with $|\gamma(n)|\le 1.$ Denote
$$
\mathbf{e}_{m}(z)= \exp(2\pi i z/m).
$$

\begin{theorem}
\label{thm:intdf} Let $a$ be any integer coprime to $m.$ For any
integers $K$ and $N$ with
$$
0<K+1\le K+N\le p-1, \quad N\ge Tp^{-1/2}(\log p)^3
$$
and any set $\mathcal{X}\subset \mathbb{Z}_{p-1},$ the inequality
$$
\sum_{x\in
\mathcal{X}}\Big|\sum_{y=K+1}^{K+N}\gamma(y)\mathbf{e}_{p}(a
\lambda^{xy})\Big|\ll\frac{|\mathcal{X}|^{1/2}N^{3/4}p^{\frac{7}{8}+o(1)}}{T^{1/4}}.
$$
holds, where $|\mathcal{X}|$ denotes the cardinality of the set
$\mathcal{X}.$
\end{theorem}

\begin{cor} \label{cor:interv} Let $a$ be any integer coprime to $m.$
For any integers $K, N, L, M$ with
$$
0<K+1\le K+N\le p-1, \quad 0< L+1\le L+M\le p-1
$$
and any coefficients $\alpha_x$ and $\beta_y$  with $|\alpha_x|\le
1,$ $|\beta_y|\le 1,$ the following inequality holds:
$$
\Big|\sum_{x=L+1}^{L+M}\sum_{y=K+1}^{K+N}\alpha_x\beta_y\mathbf{e}_{p}(a
g^{xy})\Big|\ll \left(NM\right)^{5/8}p^{5/8+o(1)}.
$$
\end{cor}
This estimate is nontrivial when $N\ge M\ge p^{5/6+\varepsilon}.$
For more information on trigonometric sums with an exponential
function and their applications, see~\cite{Ba}-\cite{Gar1} and
therein references.

To prove Corollary~\ref{cor:interv} we observe that the statement is
trivial if $N\le p^{2/3}$ or if $M\le p^{2/3}.$ Assuming $ \min\{N,
M\}>p^{2/3},$ from Theorem~\ref{thm:intdf} we derive,
$$
\Big|\sum_{x=L+1}^{L+M}\sum_{y=K+1}^{K+N}\alpha_x\beta_y\mathbf{e}_{p}(a
g^{xy})\Big|\ll\sum_{x=L+1}^{L+M}\Big|\sum_{y=K+1}^{K+N}\beta_y\mathbf{e}_{p}(a
g^{xy})\Big|\ll M^{1/2}N^{3/4}p^{5/8+o(1)},
$$
$$
\Big|\sum_{x=L+1}^{L+M}\sum_{y=K+1}^{K+N}\alpha_x\beta_y\mathbf{e}_{p}(a
g^{xy})\Big|\ll
\sum_{y=K+1}^{K+N}\Big|\sum_{x=L+1}^{L+M}\alpha_x\mathbf{e}_{p}(a
g^{xy})\Big|\ll N^{1/2}M^{3/4}p^{5/8+o(1)}.
$$
The result now follows.

In passing, we remark that Theorem~\ref{thm:intdf} and
Corollary~\ref{cor:interv} remain true if $\mathbf{e}_{p}(a g^{xy})$
is replaced by $\chi(g^{xy}+a),$ where $\chi$ is any nonprincipal
character modulo $p.$

Theorem~\ref{thm:main} also finds its application in the problem of
multiplication of intervals in a residue ring modulo $m.$
\begin{cor}
\label{cor:multinter} For any fixed $\varepsilon>0$ the set
$$
\{xy\pmod m:\quad 1\le x\le m^{1/2}, \quad S+1\le y\le
S+m^{1/2}(\log m)^{2+\varepsilon}\}
$$
contains $(1+O((\log m)^{-\varepsilon}))m$ residue classes modulo
$m.$
\end{cor}

  The classical
conjecture claims that for any prime number $p$ any nonzero residue
class modulo $p$ can be represented in the form $xy\pmod p,$ where
$1\le x, y\le p^{1/2+o(1)}.$ A weaker version of this conjecture has
been stated in~\cite{Gar2}, namely, for any prime $p$ there are
$(1+o(1))p$ residue classes modulo $p$ of the form $xy\pmod p$ with
$1\le x, y\le p^{1/2+o(1)}.$ Furthermore, in~\cite{Gar2} it has been
proved that for almost all primes $p$ almost all residue classes
modulo $p$ are representable in the form $xy\pmod p$ with $1\le x,
y\le p^{1/2}(\log p)^{1,087}.$ The following consequence of
Corollary~\ref{cor:multinter} confirms the validity of the weaker
version of the classical conjecture and essentially improves one of
our results from~\cite{GaKa}.
\begin{cor}
\label{cor:multint} For any fixed $\varepsilon>0$ and any prime
number $p$ the set
$$
\{xy\pmod p:\quad 1\le x, y\le p^{1/2}(\log p)^{2+\varepsilon}\}
$$
contains $(1+o(1))p$ residue classes modulo $p.$
\end{cor}

The method of the proof of Theorem~\ref{thm:main} combined with an
argument similar to that of~\cite{Gar3} allows to improve the
exponent of the logarithmic factor in
Corollaries~\ref{cor:multinter} and~\ref{cor:multint}. More
precisely, the following statement takes place.
\begin{theorem}
\label{thm:mltintcmp} Let $\Delta=\Delta(m)\to \infty$ as $m\to
\infty.$ Then the set
$$
\{qy\pmod m:\quad 1\le q\le m^{1/2}, \quad S+1\le y\le S+\Delta
m^{1/2}\sqrt{m/\phi(m)}\log m\}
$$
contains $(1+O(\Delta^{-1}))m$ residue classes modulo $m.$
\end{theorem}
In particular we have
\begin{cor}
\label{cor:multintprime} Let $\Delta=\Delta(p)\to \infty$ as $p\to
\infty.$ Then the set
$$
\{qy\pmod p:\quad 1\le q\le p^{1/2}, \quad 1\le y\le\Delta
p^{1/2}\log p\}
$$
contains $(1+O(\Delta^{-1}))p$ residue classes modulo $p.$
\end{cor}
Since there are $O(p^{1/2}(\log p)^{-1})$ primes not exceeding
$p^{1/2},$ we see that the set
$$
\{qy: \quad 1\le x\le p^{1/2}, \quad S+1\le y\le S+\Delta
p^{1/2}\log p\}
$$
contains only $O(p\Delta)$ integers. This shows that the ranges of
variables in Theorem~\ref{thm:mltintcmp} and
Corollary~\ref{cor:multintprime} are sufficiently sharp.

We will also prove the corresponding result for the ratio of
intervals modulo a prime which improves one of the results
of~\cite{Gar2}.
\begin{theorem}
\label{thm:ratiointprime} Let $\Delta=\Delta(p)\to \infty$ as $p\to
\infty.$ Then the set
$$
\{xy^{-1}\pmod p:\quad N+1\le x \le N+\Delta p^{1/2}, \quad S+1\le
y\le S+\Delta p^{1/2}\}
$$
contains $(1+O(\Delta^{-2}))p$ residue classes modulo $p.$
\end{theorem}

Note however, that when $N=S=0$ and $\Delta<p^{1/2}/2,$ the set
described in Theorem~\ref{thm:ratiointprime} misses $\gg
p^{1/2}\Delta^{-1}$ reside classes modulo $p,$  see~\cite{Gar2}.

For the detailed description on the multiplication table problem
modulo a prime, see~\cite{Gar2} and also~\cite{GaKa}.

The proofs of the results of~\cite{Gar2} and~\cite{GaKa} are based
on estimates of multiplicative character sums. The approach we use
here is based on trigonometric sums.

\section{Proof of Theorem~\ref{thm:main}}

Recall that $J$ denotes the number of solutions to the congruence
$$
v_1y_1\equiv v_2y_2\pmod m,\quad v_1, v_2\in \mathcal{V},\quad
S+1\le y_1, y_2\le S+L.
$$
We express $J$ in terms of trigonometric sums. Since
$$
v_1v_2^{-1}y_1\equiv y_2\pmod m,
$$
then
$$
J=\frac{1}{m}\sum_{a=0}^{m-1}\sum_{v_1\in \mathcal{V}}\sum_{v_2\in
\mathcal{V}}\sum_{y_1\in I}\sum_{y_2\in I}\mathbf{e}_{m}(a
(v_1v_2^{-1}y_1-y_2)),
$$
where $I$ denotes the interval $[S+1, S+L].$ Picking up the term
corresponding to $a=0,$ we obtain
$$
J=
\frac{|\mathcal{V}|^2L^2}{m}+\frac{1}{m}\sum_{a=1}^{m-1}\sum_{v_1\in
\mathcal{V}}\sum_{v_2\in \mathcal{V}}\sum_{y_1\in I}\sum_{y_2\in
I}\mathbf{e}_{m}(a (v_1v_2^{-1}y_1-y_2)).
$$
Furthermore,
\begin{eqnarray*}
\lefteqn{ \frac{1}{m}\sum_{a=1}^{m-1}\sum_{v_1\in
\mathcal{V}}\sum_{v_2\in \mathcal{V}}\sum_{y_1\in I}\sum_{y_2\in
I}\mathbf{e}_{m}(a
(v_1v_2^{-1}y_1-y_2))=}\\
& &\frac{1}{m}\sum_{a=1}^{m-1}\sum_{v\in \mathcal{V}}\sum_{y_1\in
I}\sum_{y_2\in
I}\mathbf{e}_{m}(a(y_1-y_2))+\\
& &\frac{1}{m}\sum_{a=1}^{m-1}\sum_{v_1\in
\mathcal{V}}\sum_{\substack{v_2\in
\mathcal{V}\\v_2\not=v_1}}\sum_{y_1\in I}\sum_{y_2\in
I}\mathbf{e}_{m}(a (v_1v_2^{-1}y_1-y_2))=\\
& &
|\mathcal{V}|L-\frac{|\mathcal{V}|L^2}{m}+\frac{1}{m}\sum_{a=1}^{m-1}\sum_{v_1\in
\mathcal{V}}\sum_{\substack{v_2\in
\mathcal{V}\\v_2\not=v_1}}\sum_{y_1\in I}\sum_{y_2\in
I}\mathbf{e}_{m}(a (v_1v_2^{-1}y_1-y_2)).
\end{eqnarray*}
Therefore
$$
J=\frac{|\mathcal{V}|^2L^2}{m}+|\mathcal{V}|L-\frac{|\mathcal{V}|L^2}{m}+
\frac{\theta_1}{m}\sum_{a=1}^{m-1}\sum_{v_1\in
\mathcal{V}}\sum_{\substack{v_2\in
\mathcal{V}\\v_2\not=v_1}}\left|\sum_{y_1\in I}\sum_{y_2\in
I}\mathbf{e}_{m}(a (v_1v_2^{-1}y_1-y_2))\right|.
$$
Here $\theta_1,$ and $\theta_j$ everywhere below, denote some
functions with $|\theta_j|\le 1.$

For a given $n$ let $r(n)$ be the number of solutions of the
congruence
$$
v_1v_2^{-1}\equiv n\pmod m, \quad v_1, v_2\in \mathcal{V}, \quad
v_1\not=v_2.
$$
In particular $r(1)=0,$ and if $(n,m)>1,$ then $r(n)=0.$ Therefore,
the above formula takes the form
$$
J=\frac{|\mathcal{V}|^2L^2}{m}+|\mathcal{V}|L-\frac{|\mathcal{V}|L^2}{m}+
\frac{\theta_1}{m}\sum_{a=1}^{m-1}\sum_{\substack{1\le n\le m\\(
n,m)=1}}r(n)\left|\sum_{y_1\in I}\sum_{y_2\in I}\mathbf{e}_{m}(a
(ny_1-y_2))\right|.
$$
It is important to note that $v^2\le m$ for any $v\in \mathcal{V}.$
For this reason we have $r(n)\le 1$ for any $n, 1\le n\le m.$
Indeed, if
$$
v_1v_2^{-1}\equiv v_3v_4^{-1}\pmod m
$$
for some $v_1, v_2, v_3, v_4\in \mathcal{V}$ and if $v_1\not=v_2,$
then
$$
v_1v_4\equiv v_3v_2\pmod m.
$$
Since $v^2\le m$ for any $v\in \mathcal{V},$ then we derive that
$v_1v_4= v_3v_2.$ But $\mathcal{V}$ consists only on prime numbers
and $v_1\not=v_2.$ Hence, $v_1=v_3, v_2=v_4.$

Thus
\begin{equation}
\label{eqn:asympJ}
J=\frac{|\mathcal{V}|^2L^2}{m}+|\mathcal{V}|L-\frac{|\mathcal{V}|L^2}{m}+
\frac{\theta_2}{m}\sum_{a=1}^{m-1}\sum_{\substack{1\le n\le m\\(
n,m)=1}}\left|\sum_{y_1\in I}\sum_{y_2\in I}\mathbf{e}_{m}(a
(ny_1-y_2))\right|.
\end{equation}
It is now useful to recall the bound
$$
\left|\sum_{y\in I}\mathbf{e}_{m}(by)\right|\le \frac{1}{|\sin(\pi
b/m|},
$$
which, in application to~\eqref{eqn:asympJ}, yields
\begin{equation}
\label{eqn:sinus}
J=\frac{|\mathcal{V}|^2L^2}{m}+|\mathcal{V}|L-\frac{|\mathcal{V}|L^2}{m}+
\frac{\theta_3}{m}\sum_{a=1}^{m-1}\sum_{\substack{1\le n\le m\\
(n,m)=1}}\frac{1}{|\sin(\pi an/m|}\frac{1}{|\sin(\pi a/m|}.
\end{equation}
For each divisor $s|m$ we collect together the values of $a$ with
$(a,m)=s.$ Then
\begin{eqnarray*}
\lefteqn{\sum_{a=1}^{m-1}\sum_{\substack{1\le n\le m\\
(n,m)=1}}\frac{1}{|\sin(\pi an/m|}\frac{1}{|\sin(\pi a/m|}=}\\
&&\sum_{s|m}\sum_{\substack{1\le a\le m-1\\ (a, m)=s}}\sum_{\substack{1\le n\le m\\
(n,m)=1}}\frac{1}{|\sin(\pi an/m|}\frac{1}{|\sin(\pi a/m|}\le \\
&&\sum_{\substack{s|m\\s<m}}s\sum_{\substack{1\le b\le m/s-1\\ (b, m/s)=1}}\sum_{\substack{1\le n\le m/s\\
(n,m/s)=1}}\frac{1}{|\sin(\pi bn/(m/s)|}\frac{1}{|\sin(\pi
b/(m/s)|}\le\\
&&\sum_{\substack{s|m\\s<m}}s\left(\sum_{\substack{1\le b\le m/s\\
(b, m/s)=1}}\frac{1}{|\sin(\pi b/(m/s)|}\right)^2\ll\\
&&\sum_{\substack{s|m\\s<m}}s\left(\sum_{1\le b\le
m/2s}\frac{m}{bs}\right)^2\le \frac{m^{3}\log^2m}{\phi(m)},
\end{eqnarray*}
where we have used the inequality
$$
\sum_{s|m}\frac{1}{s}\le
\prod_{p|m}\frac{1}{1-p^{-1}}=\frac{m}{\phi(m)}.
$$
 Inserting this bound
into~\eqref{eqn:sinus}, we obtain the required estimate.

\section{Proof of Theorem~\ref{thm:intdf}}

If $|\mathcal{X}|^{1/2}N^{1/4}T^{1/4}p^{-7/8}\le 10$ then the
estimate of Theorem~\ref{thm:intdf} becomes trivial. Therefore, we
can suppose that
$$
Q:=|\mathcal{X}|^{1/2}N^{1/4}T^{1/4}p^{-7/8}\ge 10.
$$

For a given divisor $d\mid p-1$ let $\mathcal{L}_d$ be the set of
all integers $y$ such that $dy\in [K+1, K+N]$ and $(dy, p-1)=d.$
Then
\begin{equation}
\label{eqn:overd}
 W_{a}(\gamma; T; \mathcal{X}; K, N)\le \sum_{d\mid
p-1}\mathcal{R}_a(\gamma; d, \mathcal{X}, \mathcal{L}_d),
\end{equation}
where $ \mathcal{R}_a(\gamma; d, \mathcal{X}, \mathcal{L}_d)$ is
defined by
$$
\mathcal{R}_a(\gamma; d, \mathcal{X}, \mathcal{L}_d)=\sum_{x\in
\mathcal{X}}\Big|\sum_{y\in
\mathcal{L}_d}\gamma(dy)\mathbf{e}_{p}(ag^{tdxy})\Big|
$$
Note that $|\mathcal{L}_d|\le Nd^{-1}+1$ and that the elements of
$\mathcal{L}_d$ are relatively prime to $(p-1)/d.$

For the divisors $d\mid p-1$ with the condition $d\ge Q$ we use the
trivial estimate
$$
\mathcal{R}_a(\gamma; d, \mathcal{X}, \mathcal{L}_d)\le
|\mathcal{X}||\mathcal{L}_d|\le
\frac{|\mathcal{X}|N}{Q}+|\mathcal{X}|\ll
\frac{|\mathcal{X}|^{1/2}N^{3/4}p^{7/8}}{T^{1/4}}.
$$
Therefore, in view of $\tau(p-1)\le p^{o(1)},$
from~\eqref{eqn:overd} we obtain
\begin{equation}
\label{eqn:overdQ}
 W_{a}(\gamma; T; \mathcal{X}; K, N)\ll \sum_{\substack{d\mid
p-1\\d\le Q}}\mathcal{R}_a(\gamma; d, \mathcal{X},
\mathcal{L}_d)+\frac{|\mathcal{X}|^{1/2}N^{3/4}p^{7/8+o(1)}}{T^{1/4}}.
\end{equation}

Below we suppose that $d\le Q$ (and therefore $|\mathcal{L}_d|\le
2Nd^{-1}$) and our aim is to obtain a suitable upper bound for
$\mathcal{R}_a(\gamma; d, \mathcal{X}, \mathcal{L}_d).$

Observe that if $p^{1/4}T^{1/2}N^{-1/2}Q^{-1/2}\le \log p$ then the
estimate of Theorem~\ref{thm:intdf} becomes trivial. Indeed, in this
case we would have that
$$
p^{1/4}T^{1/2}N^{-1/2}|\mathcal{X}|^{-1/4}N^{-1/8}T^{-1/8}p^{7/16}\le
\log p.
$$
Therefore,
$$
T^{1/4}\le |\mathcal{X}|^{1/6}N^{5/12}p^{-11/24+o(1)},
$$
whence
$$
\frac{|\mathcal{X}|^{1/2}N^{3/4}p^{7/8+o(1)}}{T^{1/4}}=
|\mathcal{X}|^{1/3}N^{1/3}p^{4/3+o(1)}\ge |\mathcal{X}|N.
$$

Hence, without loss of generality we may assume that
\begin{equation}
\label{eqn:condVint} p^{1/4}T^{1/2}N^{-1/2}Q^{-1/2}\ge \log p.
\end{equation}

Denote by $\mathcal{V}$ the set of the first
$[p^{1/4}T^{1/2}N^{-1/2}d^{-1/2}]$ prime numbers which are not
divisible by $p-1.$ Since any positive integer $m$ has only $O(\log
m)$ (even $O(\log m/\log\log m)$) different prime divisors, then
from ~\eqref{eqn:condVint} we deduce that for any $v\in \mathcal{V}$
we have
$$
v\ll (|\mathcal{V}|+\log p)\log p\ll |\mathcal{V}|\log p.
$$
Here $|\mathcal{V}|,$ as before, denotes the cardinality of
$\mathcal{V},$ that is
$$
|\mathcal{V}|=[p^{1/4}T^{1/2}N^{-1/2}d^{-1/2}].
$$
We observe that if $TN\le p^{3/2,}$ then the estimate of
Theorem~\ref{thm:intdf} again becomes trivial. Hence, we may suppose
that
\begin{equation}
\label{eqn:TN} TN\ge p^{3/2}.
\end{equation}

Now we follow the idea of~\cite{Gar1} in order to relate the problem
of obtaining an upper bound for $\mathcal{R}_a(\gamma; d,
\mathcal{X}, \mathcal{L}_d)$ with Theorem~\ref{thm:main}. For a
given divisor $d\mid p-1$ denote by $\mathcal{U}_d$ the set of all
elements of the ring $\mathbb{Z}_{(p-1)/d}$ relatively prime to
$(p-1)/d,$ that is $\mathcal{U}_d=\mathbb{Z}_{(p-1)/d}^*.$ For any
given integer $y$ with $(y, (p-1)/d)=1$ consider the congruence
\begin{equation}\label{congruence}
uv\equiv y\pmod{((p-1)/d)}, \quad u\in \mathcal{U}_d,\quad v\in
\mathcal{V}.
\end{equation}
The number of solutions of this congruence is exactly equal to
$|\mathcal{V}|.$ This follows from the fact that once $v$ is fixed
then $u$ is determined uniquely.

We replace $\lambda$ by $g^{t},$ where $t=(p-1)/T,$ and consider the
sum
$$
\sum_{y\in \mathcal{L}_d}\gamma(dy)\mathbf{e}_{p}(ag^{tdxy}).
$$
Denote by $\delta(y):=\delta(\mathcal{L}_d; y)$ the characteristic
function of the set $\mathcal{L}_d$ in the ring
$\mathbb{Z}_{(p-1)/d}.$ Since the number of solutions of the
congruence~(\ref{congruence}) is equal to $|\mathcal{V}|$ for any
fixed $y\in \mathcal{L}_d,$ then
$$
\sum_{y\in
\mathcal{L}_d}\gamma(dy)\mathbf{e}_{p}(ag^{tdxy})=\frac{1}{|\mathcal{V}|}\sum_{u\in
\mathcal{U}_d}\sum_{v\in
\mathcal{V}}\gamma(duv)\delta(uv)\mathbf{e}_{p}(ag^{tdxuv}).
$$
Therefore, setting
$$
\mathcal{R}_a(\gamma; d, \mathcal{V}, \mathcal{L}_d)=\sum_{x\in
\mathcal{X}}\Big|\sum_{y\in
\mathcal{L}_d}\gamma(dy)\mathbf{e}_{p}(ag^{tdxy})\Big|
$$
we see that
$$
\mathcal{R}_a(\gamma; d, \mathcal{X},
\mathcal{L}_d)=\frac{1}{|\mathcal{V}|}\sum_{x\in
\mathcal{X}}\Big|\sum_{u\in \mathcal{U}_d}\sum_{v\in
\mathcal{V}}\gamma(duv)\delta(uv)\mathbf{e}_{p}(ag^{tdxuv})\Big|.
$$
Application of the Cauchy inequality to the sums over $x$ and $u$
yields
$$
\Big|\mathcal{R}_a(\gamma; d, \mathcal{X}, \mathcal{L}_d)\Big|^2\le
\frac{|\mathcal{U}_d||\mathcal{X}|}{|\mathcal{V}|^2} \sum_{u\in
\mathcal{U}_d}\sum_{x=1}^{p-1} \Big|\sum_{v\in\mathcal{
V}}\gamma(duv)\delta(uv) \mathbf{e}_{p}(ag^{tdxuv})\Big|^2.
$$
If $(n, p-1)=d,$ if $x$ runs through $\mathbb{Z}_{p-1}$ and if $z$
runs through the reduced residue system modulo $p,$ then $g^{nx}$
and $z^{d}$ run the same system of residues modulo $p$ (including
the multiplicities). Since $(du, p-1)=d,$ then
$$
\Big|\mathcal{R}_a(\gamma; d, \mathcal{X}, \mathcal{L}_d)\Big|^2\le
\frac{|\mathcal{U}_d||\mathcal{X}|}{|\mathcal{V}|^2} \sum_{u\in
\mathcal{U}_d}\sum_{z=1}^{p-1}\Big|\sum_{v\in
\mathcal{V}}\gamma(duv)\delta(uv)\mathbf{e}_{p}( az^{tdv})\Big|^2,
$$
whence
\begin{eqnarray*}
&&\lefteqn{\Big|\mathcal{R}_a(\gamma; d, \mathcal{X}, \mathcal{L}_d)\Big|^2 \le} \\
& &  \frac{|\mathcal{U}_d||\mathcal{X}|}{|\mathcal{V}|^2} \sum_{u\in
\mathcal{U}_d}\sum_{v_1\in \mathcal{V}}\sum_{v_2\in
\mathcal{V}}\gamma(duv_1)\overline{\gamma(duv_2)}\delta(uv_1)\delta(uv_2)\sum_{z=1}^{p-1}\mathbf{e}_{p}(
a(z^{tdv_1}-z^{tdv_2})).
\end{eqnarray*}

The rightmost sum is equal to $p-1$ when $v_1=v_2$ and, according to
the Weil estimate, is bounded by $(\max\{v_1, v_2\})tdp^{1/2}$ when
$v_1\not=v_2.$ Recall that
$$
\max\{v_1, v_2\}\ll |\mathcal{V}|\log p,
$$
and $|\mathcal{U}_d|=\phi(\frac{p-1}{d})\le \frac{p}{d}.$ Therefore,
\begin{eqnarray}
\label{eqn:RdXLd}
\begin{split}
\lefteqn{
\Big|\mathcal{R}_a(\gamma; d, \mathcal{X}, \mathcal{L}_d)\Big|^2} \\
& & \ll \frac{p^2|\mathcal{X}|}{d|\mathcal{V}|^2} \sum_{u\in
\mathcal{U}_d}\sum_{v\in \mathcal{V}}\delta^2(uv)+
\frac{p^{3/2}|\mathcal{X}|t\log p}{|\mathcal{V}|} \sum_{u\in
\mathcal{U}_d}\sum_{v_1\in \mathcal{V}}\sum_{v_2\in
\mathcal{V}}\delta(uv_1)\delta(uv_2).
\end{split}
\end{eqnarray}
Next, from~(\ref{congruence}) we derive the formula
\begin{equation}\label{eqn:delta1}
 \sum_{u\in \mathcal{U}_d}\sum_{v\in \mathcal{V}}\delta^2(uv)=|\mathcal{V}||\mathcal{L}_d|.
\end{equation}
Now set
$$
J=\sum_{u\in \mathcal{U}_d}\sum_{v_1\in \mathcal{V}}\sum_{v_2\in
\mathcal{V}}\delta(uv_1)\delta(uv_2)
$$
and observe that $J$ is equal to the number of solutions of the
system of congruences
$$\left\{
  \begin{array}{ll}
    uv_1\equiv y_2\pmod {\frac{p-1}{d}}\\
    uv_2\equiv y_1\pmod {\frac{p-1}{d}}
  \end{array}
\right.
$$
subject to the conditions
$$
u\in U_d, \quad v_1, v_2\in \mathcal{V},\quad y_1, y_2\in
\mathcal{L}_d.
$$
It then follows that
$$
v_1y_1\equiv v_2y_2\pmod {\frac{p-1}{d}}.
$$
Therefore, from~\eqref{eqn:RdXLd} and~\eqref{eqn:delta1}, we derive
that
$$
\Big|\mathcal{R}_a(\gamma; d, \mathcal{X}, \mathcal{L}_d)\Big|^2 \ll
\frac{p^2|\mathcal{X}|}{d|\mathcal{V}|^2}|\mathcal{V}||\mathcal{L}_d|
+ \frac{p^{5/2}|\mathcal{X}|\log p}{T|\mathcal{V}|}J,
$$
whence
\begin{equation}
\label{eqn:intRdXLd} \Big|\mathcal{R}_a(\gamma; d, \mathcal{X},
\mathcal{L}_d)\Big|^2 \ll \frac{p^2|\mathcal{X}|N}{d^2|\mathcal{V}|}
+ \frac{p^{5/2}|\mathcal{X}|\log p}{T|\mathcal{V}|}J_d,
\end{equation}
where $J_d$ denotes the number of solutions of the congruence
\begin{equation*}
v_1y_1\equiv v_2y_2\pmod {\frac{p-1}{d}},\quad v_1, v_2\in
\mathcal{V}, \quad \frac{K+1}{d}\le y_1, y_2\le \frac{K+N}{d}.
\end{equation*}

It is important to note that the condition of
Theorem~\ref{thm:intdf} yields, for any $v\in \mathcal{V},$ the
bound
$$
v^2\le |\mathcal{V}|^2(\log p)^{2+o(1)}\le p^{1/2}(\log
p)^{2+o(1)}TN^{-1}d^{-1}\le \frac{p-1}{d}.
$$
Hence, we can apply Theorem~\ref{thm:main} with $m=(p-1)/d.$ It
gives
$$
J_d\ll \frac{|\mathcal{V}|^2N^2}{(p-1)d}+\frac{|\mathcal{V}|N}{d}+
\frac{|\mathcal{V}|N^2}{pd}+\frac{p(\log\log p)\log^2p}{d},
$$
whence, using  $|\mathcal{V}|\le p^{1/4}T^{1/2}N^{-1/2},$ we obtain
$$
J_d\ll
\frac{p^{-1/2+o(1)}TN}{d}\left(1+\frac{p^{3/4}}{T^{1/2}N^{1/2}}+\frac{p^{3/2}}{TN}\right).
$$
Taking into account~\eqref{eqn:TN}, we deduce
$$
J_d\ll \frac{p^{-1/2+o(1)}TN}{d}.
$$
Combining this estimate with~\eqref{eqn:intRdXLd}, we conclude
$$
\Big|\mathcal{R}_a(\gamma; d, \mathcal{X},
\mathcal{L}_d)\Big|\ll\frac{|\mathcal{V}|^{1/2}N^{3/4}p^{7/8+o(1)}}{d^{1/4}T^{1/4}}.
$$
The result now follows in view of~\eqref{eqn:overdQ}.

\section{Proof of Theorem~\ref{thm:mltintcmp}}

Without loss of generality we may suppose that
$$
\Delta m^{1/2}\sqrt{m/\phi(m)}\log m<m,
$$
as otherwise the statement of Theorem~\ref{thm:mltintcmp} is
trivial.

Denote by $\mathcal{V}$ the set of prime numbers coprime to $m$ and
not exceeding $m^{1/2}.$ Let $J$ denote the number of solutions to
the congruence
$$
v_1(y_1+z_1)\equiv v_2(y_2+z_2)\pmod m
$$
subject to the conditions
$$
\quad v_1, v_2\in \mathcal{V},\quad y_1, y_2, z_1, z_2\in I,
$$
where $I$ denotes the set of integers $x, [S/2]+1\le x\le [S/2]+L,$
and
$$
L=\left[\frac{\Delta m^{1/2}\sqrt{m/\phi(m)}\log m}{2}\right].
$$

Obviously that
$$
S+1\le y_i+z_i\le S+\Delta m^{1/2}\log m, \quad i=1, 2.
$$
Following the lines of the proof of Theorem~\ref{thm:main}, we
express $J$ in terms of trigonometric sums. Since
$$
v_1v_2^{-1}(y_1+z_1)\equiv y_2+z_2\pmod m,
$$
then
$$
J=\frac{1}{m}\sum_{a=0}^{m-1}\sum_{v_1\in \mathcal{V}}\sum_{v_2\in
\mathcal{V}}\sum_{y_1, z_1\in I}\sum_{y_2, z_2\in I}\mathbf{e}_{m}(a
(v_1v_2^{-1}(y_1+z_1)-y_2-z_2)).
$$
Picking up the term corresponding to $a=0,$ we obtain
$$
J=
\frac{|\mathcal{V}|^2L^4}{m}+\frac{1}{m}\sum_{a=1}^{m-1}\sum_{v_1\in
\mathcal{V}}\sum_{v_2\in \mathcal{V}}\sum_{y_1, z_1\in I}\sum_{y_2,
z_2\in I}\mathbf{e}_{m}(a (v_1v_2^{-1}(y_1+z_1)-y_2-z_2)).
$$
Since the number of solutions of the congruence
$$
y_1+z_1\equiv y_2+z_2\pmod m, \quad y_1, z_1, y_2, z_2 \in I
$$
is $O(L^3),$ then we obtain
\begin{eqnarray*}
\lefteqn{\frac{1}{m}\left|\sum_{a=1}^{m-1}\sum_{v\in
\mathcal{V}}\sum_{y_1, z_1\in I}\sum_{y_2, z_2\in
I}\mathbf{e}_{m}(a(y_1+z_1-y_2-z_2))\right|\le}\\
&&\frac{|\mathcal{V}|}{m}\sum_{a=0}^{m-1}\left|\sum_{y\in
I}\mathbf{e}_{m}(ay_1)\right|^4\ll |\mathcal{V}|L^3.
\end{eqnarray*}
Therefore,
\begin{eqnarray*}
\frac{1}{m}\sum_{a=1}^{m-1}\sum_{v_1\in \mathcal{V}}\sum_{v_2\in
\mathcal{V}}\sum_{y_1, z_1\in I}\sum_{y_2, z_2\in I}\mathbf{e}_{m}(a
(v_1v_2^{-1}(y_1+z_1)-y_2-z_2))=\\
O(|\mathcal{V}|L^3)+\frac{1}{m}\sum_{a=1}^{m-1}\sum_{v_1\in
\mathcal{V}}\sum_{\substack{v_2\in
\mathcal{V}\\v_2\not=v_1}}\sum_{y_1\in I}\sum_{y_2\in
I}\mathbf{e}_{m}(a (v_1v_2^{-1}(y_1+z_1)-y_2-z_2)).
\end{eqnarray*}
Using exactly the same argument that we used in the proof of
Theorem~\ref{thm:main}, we derive the formula
$$
J=\frac{|\mathcal{V}|^2L^4}{m}+ O(|\mathcal{V}|L^3)+O(R),
$$
where
$$
R=\frac{1}{m}\sum_{a=1}^{m-1}\sum_{\substack{1\le n\le m\\(
n,m)=1}}\left|\sum_{y_1, z_1\in I}\sum_{y_2, z_2\in
I}\mathbf{e}_{m}(a (n(y_1+z_1)-y_2-z_2))\right|
$$
Next, introducing $s=(a,m),$ we obtain
\begin{eqnarray*}
R=\frac{1}{m}\sum_{\substack {s|m\\s<m}}\sum_{\substack {b\le m/s-1\\
(b, m/s)=1}}\sum_{\substack{1\le n\le m\\( n,m)=1}}\left|\sum_{y_1,
z_1\in I}\sum_{y_2, z_2\in I}\mathbf{e}_{m/s}(b
(n(y_1+z_1)-y_2-z_2))\right|\le \\
\frac{1}{m}\sum_{\substack {s|m\\s<m}}s\sum_{\substack {b\le m/s-1\\
(b, m/s)=1}}\sum_{\substack{1\le n\le m/s\\(
n,m/s)=1}}\left|\sum_{y_1, z_1\in I}\sum_{y_2, z_2\in
I}\mathbf{e}_{m/s}(bn(y_1+z_1)-b(y_2+z_2))\right|\le\\
\frac{1}{m}\sum_{\substack {s|m\\s<m}}s\left(\sum_{\substack{1\le
n\le m/s\\( n,m/s)=1}}\left|\sum_{y_1, z_1\in
I}\mathbf{e}_{m/s}(n(y_1+z_1))\right|\right)^2=\\
\frac{1}{m}\sum_{\substack {s|m\\s<m}}s\left(\sum_{\substack{1\le
n\le m/s\\( n,m/s)=1}}\left|\sum_{y\in
I}\mathbf{e}_{m/s}(ny)\right|^2\right)^2.
\end{eqnarray*}
Therefore,
\begin{equation}
\label{eqn:Jcmp} J=\frac{|\mathcal{V}|^2L^4}{m}+
O(|\mathcal{V}|L^3)+O(R_1)+O(R_2),
\end{equation}
where
\begin{equation}
\label{eqn:R1int} R_1=\frac{1}{m}\sum_{\substack
{s|m\\s<m/L}}s\left(\sum_{\substack{1\le n\le m/s\\(
n,m/s)=1}}\left|\sum_{y\in I}\mathbf{e}_{m/s}(ny)\right|^2\right)^2
\end{equation}
\begin{equation}
\label{eqn:R2int} R_2=\frac{1}{m}\sum_{\substack {s|m\\m/L\le
s<m}}s\left(\sum_{\substack{1\le n\le m/s\\(
n,m/s)=1}}\left|\sum_{y\in I}\mathbf{e}_{m/s}(ny)\right|^2\right)^2
\end{equation}
If $s<m/L$ then $m/s>L$ and therefore, the congruence
$$
y_1\equiv y_2\pmod {m/s}, \quad y_1, y_2\in I
$$
has $L$ solutions. Hence,
$$
\sum_{1\le n\le m/s}\left|\sum_{y\in
I}\mathbf{e}_{m/s}(ny)\right|^2=\frac{mL}{s},
$$
whence, using~\eqref{eqn:R1int},
\begin{eqnarray*}
 R_1\le \frac{1}{m}\sum_{\substack
{s|m\\s<m/L}}s\left(\sum_{1\le n\le m/s}\left|\sum_{y\in
I}\mathbf{e}_{m/s}(ny)\right|^2\right)^2=\\
mL^2\sum_{\substack {s|m\\s<m/L}}s^{-1}\le mL^2\sum_{s|m}s^{-1}\le
\frac{m^2L^2}{\phi(m)}.
\end{eqnarray*}
Inserting this bound into~\eqref{eqn:Jcmp}, we deduce
\begin{equation}
\label{eqn:Jcmp2} J=\frac{|\mathcal{V}|^2L^4}{m}+
O(|\mathcal{V}|L^3)+O(m^2L^2/\phi(m))+O(R_2).
\end{equation}

Now we proceed to estimate $R_2.$ Note that in~\eqref{eqn:R2int} we
have $(n,m/s)=1.$ Therefore, for any integer $K,$
$$
\sum_{y=K+1}^{K+m/s}\mathbf{e}_{m/s}(ny)=0,
$$
whence we deduce that there exist integers $A$ and $B$ with $0< B\le
m/s$ such that
$$
\sum_{y\in I}\mathbf{e}_{m/s}(ny)=\sum_{A<y\le
A+B}\mathbf{e}_{m/s}(ny).
$$
Hence
\begin{eqnarray*}
\sum_{\substack{1\le n\le m/s\\ (n,m/s)=1}}\left|\sum_{y\in
I}\mathbf{e}_{m/s}(ny)\right|^2=\sum_{\substack{1\le n\le m/s\\
(n,m/s)=1}}\left|\sum_{A<y\le A+B}\mathbf{e}_{m/s}(ny)\right|^2\le \\
\sum_{n=1}^{m/s}\left|\sum_{A<y\le
A+B}\mathbf{e}_{m/s}(ny)\right|^2=mB/s\le m^2/s^2.
\end{eqnarray*}
Taking this into account, from~\eqref{eqn:R2int} we deduce
$$
R_2\le \frac{1}{m}\sum_{s\ge m/L}s(m^4/s^4)\ll mL^2.
$$
Therefore, in view of~\eqref{eqn:Jcmp2}, we obtain the asymptotic
formula
\begin{eqnarray*}
J=\frac{|\mathcal{V}|^2L^4}{m}+
O(|\mathcal{V}|L^3)+O(m^2L^2/\phi(m))=\\
\frac{|\mathcal{V}|^2L^4}{m}\left(1+O\left(\frac{m}{|\mathcal{V}|L}+
\frac{m^3}{\phi(m)|\mathcal{V}|^2L^2}\right)\right).
\end{eqnarray*}
Recalling that $|\mathcal{V}|\gg m^{1/2}/\log m$ and
$L=\left[\frac{\Delta m^{1/2}\sqrt{m/\phi(m)}\log m}{2}\right],$ we
arrive at the formula
$$
J=\frac{|\mathcal{V}|^2L^4}{m}\left(1+O(\Delta^{-1})\right).
$$

Next, define
$$
H=\{q(y+z)\pmod m,\quad q\le m^{1/2}, \quad [S/2]+1\le y, z\le
[S/2]+L\}.
$$
Obviously, $S+1\le y+z\le S+\Delta m^{1/2}\sqrt{m/\phi(m)}\log m.$
For a given $h\in H,$ by $J(h)$ we denote the number of solutions of
the congruence
$$
q(y+z)\equiv h\pmod m,\quad q\le m^{1/2}, \quad [S/2]+1\le y, z\le
[S/2]+L.
$$
Then
$$
J=\sum_{h\in H}J^2(h)\ge \frac{1}{|H|}\left(\sum_{h\in
H}J(h)\right)^{2}=\frac{1}{|H|}|\mathcal{V}|^2L^4.
$$
Therefore,
$$
|H|\ge\frac{|\mathcal{V}|^2L^4}{J}\ge\frac{m}{1+O(\Delta^{-1})}=(1+O(\Delta^{-1}))m.
$$
The result now follows in view of $|H|\le m.$

\section{Proof of Theorem~\ref{thm:ratiointprime}}

Without loss of generality we may suppose that
$$
0<N< N+\Delta p^{1/2}<p, \quad 0< M<M+\Delta p^{1/2}<p.
$$
Denote $X=[\Delta p^{1/2}/2],$ $N_1=[N/2],$ $S_1=[S/2],$ and let
$H^*$ be the set of all residue classes of the form
$(x+t)(y+z)^{-1}\pmod p,$ where
$$
N_1+1\le x, t\le N_1+X, \quad S_1+1\le y, z\le S_1 +X.
$$
Obviously,
$$
N+1\le x+t\le N+\Delta p^{1/2}, \quad S+1\le y+z\le S+\Delta
p^{1/2}.
$$
Next, let
$$
H_1^*=\{h\pmod p: \quad h\not\in H^*,\quad h\not\equiv 0\pmod p\}.
$$
Then the congruence
$$
x+t-(y+z)h\equiv 0\pmod p
$$
has no solutions in variables $h, x, t, y, z$ subject to the
condition
$$
h\in H_1^*, \quad N_1+1\le x, t\le N_1+X, \quad S_1+1\le y, z\le
S_1+X.
$$
Therefore,
$$
\sum_{a=0}^{p-1}\sum_{h\in H_1^*}\sum_{x, t\in I_1}\sum_{y, z\in
I_2}e^{2\pi i\frac{a(x+t-h(y+z))}{p}}=0,
$$
where $I_1$ and $I_2$ denote the intervals $[N_1+1, N_1+X]$ and
$[S_1+1, S_1+X]$ correspondingly.

 Separating the term corresponding to $a=0$ we deduce that
$$
|H_1^*|X^4\le \sum_{a=1}^{p-1}\left|\sum_{x, t\in I_1}e^{2\pi
i\frac{a(x+t)}{p}}\right|\left|\sum_{y, z\in I_2}\sum_{h\in
H_1^*}e^{2\pi i\frac{ah(y+z)}{p}}\right|.
$$
On the other hand for $(a, p)=1$ we have
\begin{eqnarray*}
\left|\sum_{y, z\in I_2}\sum_{h\in H_1^*}e^{2\pi
i\frac{ah(y+z)}{p}}\right|\le \sum_{h\in H_1^*}\left|\sum_{y, z\in
I_2}e^{2\pi i\frac{ah(y+z)}{p}}\right|\le
\\ \sum_{h=1}^{p-1}\left|\sum_{y, z\in I_2}e^{2\pi
i\frac{ah(y+z)}{p}}\right|\le\sum_{h=0}^{p-1}\left|\sum_{y, z\in
I_2}e^{2\pi i\frac{h(y+z)}{p}}\right|=pX,
\end{eqnarray*}
and similarly,
$$
\sum_{a=1}^{p-1}\left|\sum_{x, t\in I_1}e^{2\pi
i\frac{a(x+t)}{p}}\right|\le pX.
$$
Hence
$$
|H_1^*|X^4\le p^2X^2,
$$
whence
$$
|H_1^*|\le \frac{p^2}{X^2}\ll p\Delta^{-2}.
$$
Since $|H|=p-1-|H_1^*|,$ then the result follows.

\end{document}